\documentclass[12pt]{amsart}
\usepackage{amssymb}

\begin{document}
\bibliographystyle{alpha}
%\vsize=7.85in
%\hsize=14cm
%\oddsidemargin=-1.2cm\oddsidemargin=-1cm
%\evensidemargin=-1.2cm\evensidemargin=-1cm
%\hoffset 5.5truemm
%\setlength{\textwidth}{25cm}
\numberwithin{equation}{section}

\def\Label#1{\label{#1}}

\def\1#1{\ov{#1}}
\def\2#1{\widetilde{#1}}
\def\3#1{\mathcal{#1}}
\def\4#1{\widehat{#1}}

\def\s{s}
\def\k{\kappa}
\def\ov{\overline}
\def\span{\text{\rm span}}
\def\tr{\text{\rm tr}}
\def\GL{{\sf GL}}
\def\xo {{x_0}}
\def\Rk{\text{\rm Rk\,}}
\def\sg{\sigma}

\def \hn {holomorphically nondegenerate}
\def\hyp{hypersurface}
\def\prt#1{{\partial \over\partial #1}}
\def\det{{\text{\rm det}}}
\def\wob{{w\over B(z)}}
\def\co{\chi_1}
\def\po{p_0}
\def\fb {\bar f}
\def\gb {\bar g}
\def\Fb {\ov F}
\def\Gb {\ov G}
\def\Hb {\ov H}
\def\zb {\bar z}
\def\wb {\bar w}
\def \qb {\bar Q}
\def \t {\tau}
\def\z{\chi}
\def\w{\tau}
\def\Z{\zeta}

\def \T {\theta}
\def \Th {\Theta}
\def \L {\Lambda}
\def\b{\beta}
\def\a{\alpha}
\def\o{\omega}
\def\l{\lambda}

\def \im{\text{\rm Im }}
\def \re{\text{\rm Re }}
\def \Char{\text{\rm Char }}
\def \supp{\text{\rm supp }}
\def \codim{\text{\rm codim }}
\def \Ht{\text{\rm ht }}
\def \Dt{\text{\rm dt }}
\def \hO{\widehat{\mathcal O}}
\def \cl{\text{\rm cl }}
\def \bR{\mathbb R}
\def \bC{\mathbb C}
\def \bP{\mathbb P}
\def \C{\mathbb C}
\def \bL{\mathbb L}
\def \bZ{\mathbb Z}
\def \bN{\mathbb N}
\def \scrF{\mathcal F}
\def \scrK{\mathcal K}
\def \scrM{\mathcal M}
\def \cR{\mathcal R}
\def \scrJ{\mathcal J}
\def \scrA{\mathcal A}
\def \scrO{\mathcal O}
\def \scrV{\mathcal V}
\def \scrL{\mathcal L}
\def \scrE{\mathcal E}
\def \hol{\text{\rm hol}}
\def \aut{\text{\rm aut}}
\def \Aut{\text{\rm Aut}}
\def \J{\text{\rm Jac}}
\def\jet#1#2{J^{#1}_{#2}}
\def\gp#1{G^{#1}}
\def\gpo{\gp {2k_0}_0}
\def\emmp {\scrF(M,p;M',p')}
\def\rk{\text{\rm rk}}
\def\Orb{\text{\rm Orb\,}}
\def\Exp{\text{\rm Exp\,}}
\def\ess{\text{\rm Ess\,}}
\def\mult{\text{\rm mult\,}}
\def\Jac{\text{\rm Jac\,}}
\def\Span{\text{\rm span\,}}
\def\d{\partial}
\def\D{\3J}
\def\pr{{\rm pr}}
\def\dbl{[\![}
\def\dbr{]\!]}
\def\nl{|\!|}
\def\nr{|\!|}

\def \depth{\text{\rm depth\,}}
\def \D{\text{\rm Der}\,}
\def \Rk{\text{\rm Rk}\,}
\def \ima{\text{\rm im}\,}
\def \vfi{\varphi}

\title {Images of real submanifolds under finite
holomorphic mappings}

\author[P. Ebenfelt and L. P. Rothschild]{Peter Ebenfelt
and Linda P. Rothschild }\thanks{The first
author is supported in part by the NSF grant DMS-0401215.
The second author is supported in part by the NSF grant
DMS-0400880.}

\address{ Department of Mathematics, University of
California at San Diego, La Jolla, CA 92093-0112, USA}
\email{pebenfel@math.ucsd.edu, lrothschild@ucsd.edu }

%\author[Peter Ebenfelt]{Peter
%  Ebenfelt$^0$}\footnotetext{{\rm Royal Swedish
%      Academy of Sciences Research Fellow supported by a
%grant
%     from the Knut and
%   Alice Wallenberg Foundation.\newline}}
%\address{Department of Mathematics, 0112, University
%of California at San Diego, La Jolla, CA 92093-0112}
%\email{
%pebenfel@math.ucsd.edu}
%\keywords \endkeywords
%\subjclass{32H40, 32V10}
%\date{\number\year-\number\month-\number\day}
%\loadeufm

\abstract We give some results concerning the smoothness of the image of a real-analytic submanifold in complex space under the action of a finite holomorphic mapping. For instance, if the submanifold is not contained in a proper complex subvariety, we give a necessary and sufficient condition guaranteeing that its image is smooth and the mapping is transversal to the image.
\endabstract

\newtheorem{Thm}{Theorem}[section]
\newtheorem{Def}[Thm]{Definition}
\newtheorem{Cor}[Thm]{Corollary}
\newtheorem{Pro}[Thm]{Proposition}
\newtheorem{Lem}[Thm]{Lemma}
\newtheorem{Rem}[Thm]{Remark}
\newtheorem{Ex}[Thm]{Example}

\maketitle
\section{Introduction and Main Results} In this paper, we study finite holomorphic mappings of real-analytic submanifolds in
$\bC^N$.  Recall that a
germ of a holomorphic mapping  $H\colon
(\bC^N,p_0)\to(\bC^N,\tilde p_0)$ is {\em finite} at $p_0$
if $H^{-1}(\tilde p_0)\cap U=\{p_0\}$ for a sufficiently
small open neighborhood $U$ of $p_0$ in $\C^N$.  If $V
\subset
\C^N$ is a germ at $p_0$ of a real-analytic subvariety and
$H$ is a finite holomorphic mapping, then its image,
$H(V)$, is contained in a germ at $\tilde p_0$ of a
real-analytic subvariety of the same dimension. We
consider here the case where $V$ is a (germ at $p_0$ of a)
real-analytic submanifold and ask for geometric
conditions  guaranteeing that the image $\tilde V:=H(V)$ is again a
(germ at $\tilde p_0$ of a) submanifold and $H$ is transversal to $\tilde V$ at $p_0$.  Our main result
(Theorem \ref{t:image}) generalizes to higher codimension
earlier work of Baouendi and the second author (see
\cite{brimages}).  This study is partly motivated by the
recent interest in the structure of nondegenerate mappings
(e.g.\ finite holomorphic mappings) taking one
real-analytic submanifold in
$\bC^N$ into another.  We mention here only the papers
\cite{er1}, \cite{lm2}, \cite{kz1},  \cite{elz},
\cite{mmzapprox}, \cite{mmzanal}, and refer the reader to
these papers for precise results and a more extensive
bibliography.

Before stating our main result, we must first introduce
some notation. Let $M$ be a real-analytic submanifold of
codimension $d$ in $\bC^N$ with $p_0\in M$.  We let $\3 M$
be the usual complexification of $M$ in some neighborhood
of $(p_0,\bar p_0)$ in
$\bC^N\times\bC^N$; i.e.\ $\mathcal M$ is defined near
$(p_0,\bar p_0)$ in
$\bC^N\times
\bC^N$ by
$\rho_1(Z,\zeta)=\ldots=\rho_d(Z,\zeta)=0$ if $M$ is
defined near $p_0$ by
\begin{equation}\Label{e:defeq}
\rho_1(Z,\bar Z)=\ldots=\rho_d(Z,\bar Z)=0.\end{equation}
We shall also associate to a holomorphic mapping $H\colon
(\bC^{N},p_0)\to (\bC^{N},\tilde p_0)$ its
complexification $\mathcal
H\colon(\bC^N\times\bC^N,(p_0,\bar p_0))
\to(\bC^N\times\bC^N,(\tilde p_0,\bar {\tilde p}_0))$
defined by $\mathcal H(Z,\zeta)=(H(Z),\bar H(\zeta))$,
where $\bar H(\zeta):=\overline{H(\bar\zeta)}$. Observe
that the mapping $H$ sends $M$ into another real-analytic
submanifold
$\widetilde M$ if and only if the complexified mapping
$\mathcal H$ sends
$\mathcal M$ to $\widetilde{\mathcal M}$, where
$\widetilde{\mathcal M}$ is the complexification of
$\widetilde M$. Also, observe that the mapping
$H$ is finite if and only if
$\mathcal H$ is finite.  It is easy to check that a necessary condition for $H(M)$ to be smooth is that $\mathcal H(\mathcal M)$ is smooth, but the converse is not true in general. (See Remark \ref{r:31}.)

A real-analytic submanifold $M$ is called {\it generic} if
$T_pM+J(T_p M)=T_p\bC^N$ for every $p\in M$. Here, $T_p Y$
denotes the (real) tangent space at $p$ to a manifold $Y$,
and $J\colon T\bC^N\to T\bC^N$ is the complex structure on
$\bC^N$. An equivalent definition can be given in terms of
local defining equations
\eqref{e:defeq} for $M$ near $p_0$, namely
$\partial\rho_1\wedge\ldots\wedge\partial\rho_d\neq 0$ on
$M$. A generic submanifold $M$ is said to be of {\it
finite type} at $p_0$ (in the sense of Kohn and
Bloom-Graham) if the (complex) Lie algebra
$\frak g_M$ generated by all smooth $(1,0)$ and $(0,1)$
vector fields tangent to $M$ satisfies $\frak g_M(p_0)=\bC
T_{p_0}M$, where $\C T_{p_0} M$ is the complexified
tangent space to $M$.   Recall that a germ of a smooth
mapping $g\colon (\bR^k,x)\to (\bR^\ell,y)$ is said to be
{\it transversal} to a smooth submanifold $Y \subset
\bR^\ell$ at $y$ if
\begin{equation}\Label{e:trans} T_yY
+dg(T_x(\bR^k))=T_y(\bR^\ell).
\end{equation} We shall say that a holomorphic mapping
$H\colon (\bC^{N},p_0)\to (\bC^{N},\tilde p_0)$ is
transversal to a real-analytic submanifold $\widetilde
M\subset \bC^N$  at $\tilde p_0$ if it is transversal to
$M$ at $0$ as a real mapping $H\colon (\bR^{2N},p_0)\to
(\bR^{2N},\tilde p_0)$. Finally, the holomorphic mapping
$H$ is said to be {\it CR transversal} to a generic
submanifold
$\widetilde M$  at $\tilde p_0$ if
\begin{equation}\Label{e:CRtransdef} T^{1,0}_{\tilde
p_0}\widetilde M+dH(T^{1,0}_{p_0}\bC^N)=T^{1,0}_{\tilde
p_0}\bC^N.
\end{equation} 
Here $T^{1,0}_{\tilde p_0}\widetilde M$ denotes the the smooth $(1,0)$ vectors in $T_{\tilde
p_0}\bC^N$ that are tangent to $\2M$ at ${\tilde
p_0}$.

The following theorem is the main result of
this paper.

\begin{Thm}\Label{t:image} Let $M$ be  a (germ of a) real-analytic submanifold through
$p_0$ in $\bC^N$  and $H\colon (\bC^N,p_0)\to
(\bC^N,\tilde p_0)$ a germ of a finite holomorphic
mapping. Consider the two properties:

{\rm (i)} The image
$H(M)$ is a germ at $\tilde p_0$ of a real-analytic
submanifold.
\smallskip

{\rm (ii)} The complexified germ $\mathcal H$ satisfies
$\mathcal H^{-1}(\mathcal H(\mathcal M))=\mathcal M$,
where $\mathcal M$ denotes the complexification of $M$.
\smallskip

If $M$ is not contained in any proper complex analytic
subvariety through $p_0$, then
\begin{equation}\Label{e:equiv1}
\text{{\rm (i)} with $H$ is transversal to $H(M)$ at
$\tilde p_0$.} \iff \text{{\rm (ii)}}
\end{equation}

If $M$ is generic and of finite type at $p_0$, then
\begin{equation}\Label {e:equiv2}
\text{{\rm (i) with }}H(M)\  \text{\rm generic} \iff \text{{\rm (ii)} }.
\end{equation} Moreover, in the latter case,  if either
{\rm (i)} or {\rm (ii)} is satisfied, then the image
$H(M)$ is  of finite type at $\tilde p_0$, and
$H$ is CR transversal to $H(M)$ at $\tilde p_0$.
\end{Thm}

\begin{Rem}  {\rm It is well known that  if $\2 M$ is a
real-analytic submanifold and
$H$ is transversal to $\2 M$ at $\tilde p_0$, then
$H^{-1}(\2 M)$ is necessarily a real-analytic
submanifold.  Moreover, if $\2 M$ is generic and $H$ is
CR  transversal to $\2 M$, then $H^{-1}(\2 M)$ is
generic.  Theorem \ref{t:image} can be viewed as
providing  partial converses to these statements.  }
\end{Rem}

Since a smooth real hypersurface in $\C^N$ is necessarily a generic submanifold, we have the following corollary, which shows that condition (iii) of Theorem 1 in \cite{brimages} is extraneous.

\begin{Cor} \Label {c:hyper} If $M$ is a real-analytic
hypersurface of finite type at
$p_0$ in $\bC^N$  and $H\colon (\bC^N,p_0)\to
(\bC^N,\tilde p_0)$ a germ of a finite holomorphic
mapping, then $H(M)$ is a real-analytic, real submanifold if and only if  $\mathcal H^{-1}(\mathcal H(\mathcal M))=\mathcal M$,
where $\mathcal M$ denotes the complexification of $M$.
\end{Cor}

\begin{Rem} {\rm  By using Theorem 4 in \cite{brimages}, we may replace (i)
by the condition
\smallskip
(i$'$)  {\it  The image $H(M)$ is a germ at $p_0$ of a
smooth  submanifold.} }
\end{Rem}

\begin{Rem} {\rm If  (ii) in Theorem
\ref{t:image} is satisfied, then
$H^{-1}(H(M)) = M$.  However, even in the case of a
hypersurface, the latter condition does not imply (i) or
(ii) (see Remark 1.11 in \cite{brimages} for an example).
}
\end{Rem}

The following example shows that the image of a generic manifold of finite type under a finite holomorphic mapping may not be a CR manifold at $0$. (Recall that $M\subset \C^N$ is {\it CR at $0$} if the mapping $p\mapsto \dim T^{0,1}M$ is constant for $p$ in a neighborhood of $0$.  A generic manifold through $0$  is necessarily CR at $0$.)  This example therefore shows that the condition that $H(M)$ is generic cannot be omitted in \eqref{e:equiv2}.

\begin {Ex}\Label{fdex}
{\rm  Let $M\subset \C^3$ be the generic hypersurface given by
\begin{equation} M: = \{(z,w_1,w_2) \in \C^3: \im w_1 = |z|^2/2,\  \im w_2 = |z|^4/2\}
\end{equation}
and $H=(F_1,F_2,G) : (\C^3,0) \to (\C^3,0)$ be the finite mapping given by
\begin{equation}
F_1(z,w) = z, \ F_2(z,w) = w_1 + i w_2, \ G(z,w) = (w_1 - i w_2)^2
\end{equation}
Let $M\subset \C^3$ be the real submanifold given by
\begin{equation} \2M: = \{(\2z_1,\2z_2,\2w) \in \C^3:  \2w = (\1{\2z}_2+i|\2z_1|^2 +|\2z_1|^4)^2\} .
\end{equation}
It is easily checked that $\2M$ is not CR at $0$.  One can check by direct calculation that $H(M) \subset \2M$.  To see  that $H$ maps $M$ onto $\2M$,  let
$(\2z_1^0,\2z_2^0,\2w^0) \in \2M$.  Taking
$ z^0: = \2z_1^0$, $ \re w_1^0 = \re\2z_2^0 + |\2z_1^0|^4/2$,  $\re w_2^0 = \im \2z_2^0 - |\2z_1^0|^2/2,$
we have $ F_1(z^0,w^0) = \2z_1^0$ and $ F_1(z^0,w^0) = \2z_2^0$, which proves the desired surjectivity.}
\end{Ex}

Note that $\C^3$ is the lowest dimensional complex space in which one can find an example of the above type.  Indeed, for a generic submanifold in $\C^2$ to be of finite type at a point, it must be a real hypersurface, so this case is covered by Corollary \ref{c:hyper}. However, a totally real generic submanifold in $\C^2$ can be mapped onto a nongeneric submanifold in  $\C^2$, as is shown by the following example.

\begin{Ex}{\rm   Let 
$$M: = \{(w_1,w_2)\in \C^2: \im w_1 = \im w_2 = 0\}$$
and let $H(w_1,w_2): = (w_1+i w_2, (w_1-i w_2)^2)$.  Then $H$ is finite, and $H$ maps $M$ onto the surface
$$\2M: = \{\2z,\2w)\in \C^2: w= \1z^2\},$$
which again is not CR at $0$ and hence not generic.
}
\end{Ex}
\begin{Ex}\Label{r:counter} {\rm  If $M$ is contained in a
complex analytic subvariety, then the implication $\Leftarrow$
in
\eqref{e:equiv1} does not hold in general. Consider $M\subset
\bC^2$ given by
$z_2=0$ and the mapping $H(z_1,z_2)=(z_1,z_2^2)$. Observe
that the complexification $\mathcal M\subset
\bC^2\times\bC^2$ of
$M$ is the submanifold of points
$(z_1,z_2,\zeta_1,\zeta_2)$ such that $z_2=\zeta_2=0$, and
the complexified map is given by
$\mathcal H(z,\zeta)=(z_1,z_2^2,\zeta_1,\zeta_2^2)$.
Clearly, we have $\mathcal H^{-1}(\mathcal H(\mathcal
M))=\mathcal M$. The image $\tilde M:=H(M)$ is a
submanifold at $0$ (i.e.\ (i)), but $H$ is not transverse
to $\tilde M$ at $0$. However, we do not
know of any examples where (ii) holds, but (i) does not.
For further discussion about this point, see
Section \ref{s:last} of this paper.
}
\end{Ex}

As mentioned above, Theorem \ref{t:image} generalizes and
sharpens a theorem of Baouendi and the second author (see
\cite{brimages}, Theorem 1, part (B)) for the case where
$M$ is an essentially finite real-analytic hypersurface.
 We should also point out
that the conclusion in Theorem \ref{t:image} that $H$ is CR transversal to $H(M)$, when $M$ is generic and of finite type,
provided that $H(M)$ is a generic  manifold, was proved
in a recent paper \cite{er1} (see Theorem 1.1) by the
authors.
We conclude the introduction by mentioning two corollaries concerning ranks of finite 
holomorphic mappings that follow from Theorem 1.1. 

\begin{Cor}\Label {c:rk} Let $M$ be a real-analytic generic submanifold of finite type through $p_0$ in
$\bC^N$, and  $H\colon (\bC^N,p_0)\to
(\bC^N,\tilde p_0)$ a germ of a finite holomorphic
mapping.   Let $\mathcal M$ and $\mathcal H$
 denote the complexifications of $M$ and $H$, respectively. If
$\mathcal H^{-1}(\mathcal H(\mathcal M))=\mathcal M$, then
\begin{equation}
\rk \frac{\partial H}{\partial Z}(p_0)\geq \codim M,
\end{equation} where $\rk$ denotes the rank of a matrix
and $\codim M$ is the real codimension of $M$ in $\bC^N$.
\end{Cor}

By combining Theorem \ref{t:image} above with a theorem from
\cite{er1}, we obtain, as a corollary, a sufficient geometric condition for a finite
mapping to be a local biholomorphism at a given point. For this
recall that $M$ is said to be {\it finitely nondegenerate} at $p_0$
if
\begin{equation}\Label{e:fnd}
\span_\bC\left \{ L^\alpha\left(\frac{\partial
\rho^j}{\partial Z}\right)(p_0)\colon j=1,\ldots, d,\
\alpha\in
\mathbb N_+^n\right
\}=\bC^N,
\end{equation} where $\span_\bC$ denotes the vector space
spanned over $\bC$ and
$L^\alpha:=L_1^{\alpha_1}\ldots L_n^{\alpha_n}$.  Here,
$L_1,\ldots, L_n$ is a basis for the smooth
$(0,1)$ (or CR) vector fields tangent to $M$ near $p_0$,
and
 $M$ is defined locally  near $p_0$ by \eqref{e:defeq}.  A
direct consequence of Theorem \ref{t:image} above and
Theorem 1.2 in \cite{er1} is the following result.

\begin{Cor} Let $M$,  $\mathcal M$, $H$, and $\mathcal H$  be as in Corollary \ref{c:rk}.  Assume, in addition, that $M$ is  finitely nondegenerate at $p_0$.   If
$\mathcal H^{-1}(\mathcal H(\mathcal M))=\mathcal M$, then
$H$ is a local biholomorphism at $p_0$.
\end{Cor}

\section {Images of complex manifolds under finite
mappings}  The study of images of complex analytic
manifolds and varieties under finite holomorphic mappings
has a long history (see e.g.\ \cite{remmert},
\cite{gunros}, \cite{grifhar}, \cite{rudinball}).   The
proof of Theorem \ref{t:image} is mainly based on the
following result
 concerning images of complex manifolds. 

\begin{Thm}\Label{t:cplximage} Let $X$ be a complex
submanifold through $0$ in
$\bC^k$ and $f\colon (\bC^k,0)\to (\bC^k,0)$ a germ of a
finite holomorphic mapping such that
\begin{equation}\Label{e:jaccond}
\left. \det \frac{\partial f}{\partial z}\right
|_X\not\equiv 0,
\end{equation} where $z=(z_1,\ldots,z_k)$ are coordinates
in
$\bC^k$. Then, the following are equivalent:

\smallskip {\rm (a)} $f(X)$ is a germ of a manifold at $0$
and $f$ is transversal to $f(X)$.

\smallskip {\rm (b)} $f^{-1}(f(X))=X$ as germs at $0$.
\end{Thm}

\begin{Rem}  {\rm Without the assumption \eqref{e:jaccond},
condition (b) in Theorem \ref{t:cplximage} does not imply that
$f$ is transversal to $f(X)$, as is shown by the example
given in Remark \ref{r:counter} above.  We do not
know if condition (b) implies that $f(X)$ is a manifold
without assuming  \eqref{e:jaccond}.  If $X$ is of dimension one, then it is shown in Theorem \ref{t:curveimage1} below that in fact (b) does imply that $f(X)$ is a manifold even without assuming \eqref{e:jaccond}.
}
\end{Rem}

\begin{Rem} {\rm We note that without the condition of transversality in (a), condition (b) need not hold. For example, consider $X=\{(z,w)\colon w=0\}$ and the mapping $f(z,w)=(z^2+w^2,zw)$. It can be easily checked that $f(X)=X$, but $f^{-1}(f(X))\neq X$.
}
\end{Rem}

\begin{proof}[Proof of Theorem $\ref{t:cplximage}$] The
proof of (a)$\implies$(b) is immediate by the
transversality assumption (without using \eqref{e:jaccond}).
We shall prove (b)$\implies$(a). We first observe that, by the
proper mapping theorem, $f(X)$ is a complex analytic
(irreducible) subvariety, of the same dimension as $X$,
through $0$ in $\bC^k$. Let  $q$ be the
codimension of $X$ in $\bC^k$.  We choose local coordinates
$(x,y)\in \bC^p\times\bC^q$, with $p+q=k$, vanishing at the origin in $\bC^k$ such
that
$X$ is given locally by $y=0$.

 Our first claim is
that the coordinate functions $y_l$, $l=1,\ldots,q$,
belong to the ideal
$I(f(x,y))$. If $m$ denotes the multiplicity at the origin
of the finite mapping $f\colon \bC^k\to\bC^k$, then
assertion (b) implies that the $m$ points (counted with
their multiplicities) in
$f^{-1}(w)$, for an arbitrary $w\in f(X)$ sufficiently
close to the origin, are all contained in $X=\{(x,0)\}$.
Moreover, it follows from the condition \eqref{e:jaccond}
that these $m$ preimages will all be distinct
(multiplicity one) for a set of $w$ which is open and
dense in the variety $f(X)$. (The set of points $w$ in
$f(X)$ for which $f^{-1}(\{w\})$ consists of fewer than
$m$ points is  contained in the image of $\{z\colon
\det(\partial f/\partial z)(z)=0\}$, which by the
assumption \eqref{e:jaccond} does not contain $X$.) Thus,
if we let
$h\colon (\bC^p,0)\to(\bC^k,0)$ denote the holomorphic
mapping defined by $h(x)=f(x,0)$, then for an open and
dense set of
$w\in f(X)$ there are $m$ preimages of $w$ under $h$. It
follows (see e.g.\ Proposition 1 on p.\ 94 in \cite{AGV85};
see also Proposition 2.4 on p.\ 168 of \cite{GG86}) that
the multiplicity of $h$ at 0 is at least $m$, i.e.\
\begin{equation}\Label{e:multineq} m\leq \dim_\bC
\bC\{x\}/I(f(x,0)),
\end{equation} where $\bC\{x\}$ denotes the ring of
convergent power series in
$x$ (or ring of germs at $0$ of holomorphic functions  in
$\bC^p$). Consider the homomorphism $\phi\colon
\bC\{x,y\}\to\bC\{x\}$ defined by $\phi(g)(x)=g(x,0)$.
Clearly, $\phi$ is surjective and sends $I(f(x,y))$ into
$I(f(x,0))$. Hence, $\phi$ induces a surjective
homomorphism
$\phi^*\colon \bC\{x,y\}/I(f(x,y))\to
\bC\{x\}/I(f(x,0))$, so that $\dim_\bC
\bC\{x\}/I(f(x,0))\leq m$. On the other hand, $m$ is the
multiplicity of the finite mapping $f$, i.e.\
\begin{equation} m= \dim_\bC \bC\{x,y\}/I(f(x,y)),
\end{equation} and hence, by \eqref{e:multineq}, we must
have that
\begin{equation}\Label{e:multid}
\dim_\bC \bC\{x\}/I(f(x,0))=\dim_\bC \bC\{x,y\}/I(f(x,y))
\end{equation} and $\phi^*$ is an isomorphism. Since
$\phi^*(y_l)=0$, for
$l=1,\ldots,q$, we conclude that
\begin{equation}\Label{e:yinI} y_l\in I(f(x,y)), \quad
l=1,\ldots,q, \end{equation} as claimed above.

Notice that as a consequence
of \eqref{e:yinI},  $\partial f/\partial y(0)$  has rank $q$.
After a linear invertible transformation in the target space
$\bC^k$ (if necessary), we can decompose its coordinates as
$w=(\xi,\eta)\in\bC^p\times\bC^q$ and write the mapping
$f(x,y)$ as $f(x,y)=(R(x,y),S(x,y))^t$, where
$R=(R_1,\ldots,R_p)^t$,
$S=(S_1\ldots,S_q)^t$, and
\begin{equation}
\frac{\partial R}{\partial y}(0)=0_{p\times q},\quad
\frac{\partial S}{\partial y}(0)=I_{q\times q},
\end{equation} where $0_{p\times q}$ denotes the $(p\times
q)$-matrix whose entries are all 0 and $I_{q\times q}$ the
$(q\times q)$ identity matrix. Hence, we can further write
the components of the mapping as
\begin{equation}\Label{e:norm}
R(x,y)=R_0(x)+R_1(x,y)y,\quad S(x,y)=y+S_0(x)+S_1(x,y)y,
\end{equation} where $R_1(x,y)$ and $S_1(x,y)$ are
$(p\times q)$-matrix and
$(q\times q)$-matrix valued functions, respectively, with
$R_1(0)=0$ and $S_1(0)=0$. Observe that the restriction to
$y=0$ is given by $h(x)=f(x,0)=(R_0(x),S_0(x))^t$.

\begin{Lem}\Label{l:s0r0} With the notation introduced
above, the germ at $0$ of the holomorphic mapping $R_0\colon
(\bC^p,0)\to (\bC^p,0)$ is finite with multiplicity $m$
and
\begin{equation}\Label{e:S0gR0} S_0(x)=g(R_0(x)),
\end{equation} where $g\colon (\bC^p,0)\to (\bC^q,0)$ is
the germ of the holomorphic mapping  that satisfies
\begin{equation}\Label{e:g}
 g(\xi):=\frac{1}{m}\sum_{\nu=1}^m S_0(x^\nu(\xi))
\end{equation} for a  generic point $\xi\in \bC^p$; here,
$x^\nu(\xi)$ denote the $m$ distinct preimages of $\xi$
under the mapping $R_0$.
\end{Lem}

\begin{proof}[Proof of Lemma $\ref{l:s0r0}$] In view of
\eqref{e:yinI}, we have
\begin{equation}\Label{e:yinI2}
y=A(x,y)R(x,y)+B(x,y)S(x,y),
\end{equation} for some matrix valued functions $A,B$. If
we Taylor expand
$A(x,y)$ and $B(x,y)$ in $y$, writing $A(x,y)=A_0(x)+O(y)$
and
$B(x,y)=B_0(x)+O(y)$, we conclude by substituting
\eqref{e:norm} into
\eqref{e:yinI2} and setting $x=0$ that
\begin{equation}\Label{e:B0inv} B_0(0)=I_{q\times q}.
\end{equation} Similarly, by setting $y = 0$ we obtain
\begin{equation}\Label{e:S0inI}
A_0(x)R_0(x)+B_0(x)S_0(x)=0.
\end{equation} It follows from
\eqref{e:B0inv} that
$B_0(x)$ is invertible near
$0$ and therefore, by \eqref{e:S0inI}, the components of
$S_0(x)$ are in the ideal $I(R_0(x))$. In other words, we
have
$I(f(x,0))=I(R_0(x))$. Thus, by \eqref{e:multid} the
number $m$ of preimages of a generic $w\in f(X)$ under
$h\colon (\bC^p,0)\to (\bC^k,0)$ is also the multiplicity
of the mapping $R_0\colon (\bC^p,0)\to (\bC^p,0)$ as
claimed in the lemma. If we write
$w=(\xi,\eta)\in \bC^p\times\bC^q$ for a point in
$f(X)$, then for generic $\xi$ we have
\begin{equation}\Label{e:S0roots}
\eta=S_0(x^\nu(\xi)),\quad\nu=1,\ldots,m .\end{equation}
If we define $g$ by \eqref{e:g}, then $g$ extends to a
holomorphic function near $0$, since it is a symmetric
function of the roots $x^1,\ldots x^m$, and
\eqref{e:S0gR0} can be verified directly from
\eqref{e:S0roots}.
\end{proof}

We may now complete the proof of Theorem
\ref{t:cplximage}. Let $g$ be as in Lemma \ref{l:s0r0}. If we
make  the change of variables
\begin{equation}
\xi'=\xi,\quad \eta'=\eta-g(\xi)
\end{equation} in the target space $\bC^k$ then by writing
 $$g(R_0(x)+R_1(x,y)y) = g(R_0(x)) + O(|x|+|y|)y$$
and using Lemma
\ref{l:s0r0}, the mapping
$f(x,y)$ takes the form
\begin{equation}\Label{e:finalform}
f(x,y)=(R_0(x)+R_1(x,y)y, y+\tilde S_1(x,y)y)^t,
\end{equation} where $\tilde S_1(x,y)$ is a $(q\times
q)$-matrix-valued function with $\tilde S_1(0)=0$. By
\eqref{e:finalform}, the $q$-dimensional complex
subvariety $f(X)$, where
$X=\{y=0\}$, is contained in the
$q$-dimensional plane $\{\eta'=0\}$ in $\bC^k$. This proves
that
$f(X)=\{\eta'=0\}$ and hence is a submanifold at $0$. It
is obvious from the form
\eqref{e:finalform} of the mapping that $f$ is transversal
to
$f(X)$. This completes the proof of Theorem
\ref{t:cplximage}.
\end{proof}

\begin{Rem} {\rm  The condition \eqref{e:jaccond} in Theorem \ref{t:cplximage} is only used to deduce that the generic number of preimages of the mapping $f|_X\colon (X,0)\to (f(X),0)$ equals the multiplicity of the mapping $f\colon( \bC^k,0)\to (\bC^k,0)$.  In fact, these two properties are equivalent as the reader can verify.
}
\end{Rem}

\section{Proof of Theorem \ref {t:image}}

\begin{proof} Without loss of generality, we may take $p_0 =
\tilde p_0 = 0$.   We assume first that
$M$ is a real analytic submanifold that is not contained in
any proper complex subvariety of $\C^N$.
To prove the
implication $\implies$ of \eqref{e:equiv1}, we suppose that
$\widetilde M: = H(M)$ is a germ at $0$ of a real-analytic
submanifold and that
$H$ is transversal to $H(M)$ at $0$.
Observe that $\widetilde M$
is of the same dimension as
$M$, since $H$ is a finite mapping. If $\tilde
\rho(Z,\bar Z)$, where $\tilde \rho=(\tilde
\rho_1,\ldots, \tilde \rho_d)$, is a defining function for
$\widetilde M$, then the fact that $H$ is transversal implies that
$\rho(Z,\bar Z):=\tilde \rho(H(Z),\bar H(\bar Z))$ is a
defining function for $M$. By simply replacing $\bar Z$ by
$\zeta$ in the above, we conclude that
$\mathcal H^{-1}(\widetilde{\mathcal M})=\mathcal M$ which is
the assertion (ii).

 To prove the
implication $\Longleftarrow$ in \eqref{e:equiv1}, we shall
need the observation that
$$
\left. \det \frac{\partial \mathcal H}{\partial
(Z,\zeta)}\right |_{\mathcal M}\not\equiv 0.
$$ Indeed, if $\det (\left. \partial \mathcal H/\partial
(Z,\zeta))\right |_{\mathcal M}\equiv 0$, then by the specific form of $\mathcal H$ it would follow that $|\det (\left. \partial H/\partial
Z)\right |_M|^2\equiv 0$ contradicting the assumptions that $H$ is finite and $M$ is not contained in a proper complex analytic subvariety.
We may
now apply Theorem \ref{t:cplximage} with $X: = \mathcal M$ and
$f: = \mathcal H$ to conclude that $\widetilde {\mathcal M}: =
\mathcal H(\mathcal M)$ is a germ at $0$ of a manifold with
$\mathcal H$ transversal to
$\widetilde {\mathcal M}$ at $0$.  Since $\mathcal H$ is finite,
$\widetilde {\mathcal M}$  has the same dimension as
$\mathcal M$. Moreover, $\widetilde {\mathcal M}$
satisfies the ``reality" symmetry: if
$(Z,\zeta)\in \widetilde {\mathcal M}$, then
$(\bar\zeta,\bar Z)\in \widetilde {\mathcal M}$. The
latter is easily verified (and the verification is left to
the reader) from the fact that $\mathcal M$ has this
symmetry, by using the specific form of
$\mathcal H$. The symmetry implies that one can find
defining equations for
$\widetilde {\mathcal M}$ near $0$ of the form $\tilde
\rho(Z,\zeta)=0$, where
$\tilde \rho=(\tilde
\rho_1,\ldots,\tilde \rho_d)$, and for each $1\leq j\leq
d$ we have $\tilde \rho_j(Z,\bar Z)$ is real-valued. It
follows that the real-valued equation $\tilde
\rho(Z,\bar Z)=0$ defines a real-analytic submanifold
$\widetilde M\subset \bC^N$ through $0$ of codimension $d$. By
construction, $H$ sends $M$ into $\widetilde M$. We must
show that $H$ sends
$M$ onto $\widetilde M$ in the sense of germs at $0$. The
fact that $\mathcal H$ is transversal to $\widetilde
{\mathcal M}$ at $0$ means that
$\rho(Z,\zeta):=\tilde \rho(H(Z),\bar H(\zeta))$ is a
defining function for $\mathcal M$ at $0$. Clearly, this
also means that $\rho(Z,\bar Z)$ is a defining function for
$M$ at $0$ and, hence, $H^{-1}(\widetilde M)=M$ as germs
at $0$. Since any representative of the germ $H$ near $0$
is an open mapping, we conclude that
$H(M)=\widetilde M$ as germs at $0$. This completes the proof
of the implication $\Longleftarrow$ in \eqref{e:equiv1}.

To complete the proof of Theorem \ref{t:image}, assume that
$M$ is a generic submanifold of finite type at $0$. We shall
show that \eqref{e:equiv2} holds. Note first that $M$ generic
implies that
$M$ is not contained in any proper complex submanifold of
$\C^N$.  It follows from \eqref{e:equiv1} that (ii) $\implies$
(i) and $H$ is transversal to $H(M)$ at $0$. As above, we note that  $\rho(Z,\bar Z)=\tilde \rho(H(Z),\overline{H(Z)})$, where  $\tilde \rho(Z,\bar Z)$ is a defining function for $H(M)$ near $0$, is a defining function for $M$ near  $0$. Now, by the chain rule
$$
\frac{\partial \rho}{\partial Z}(0)=\frac{\partial\tilde \rho}{\partial Z}(0)\frac{\partial H}{\partial Z}(0)
$$
and, hence, the rank of $(\partial\tilde \rho/\partial Z)(0)$ must be $d$ since $M$ is generic. Consequently, $H(M)$ is generic. 

For the implication  $\implies$ in \eqref{e:equiv2}, let $\widetilde
M: = H(M)$ be a generic submanifold.  It follows from Proposition 2.3 in \cite{er1}  that $\widetilde M$ is of finite type at $0$, since $H$ is
finite.  Finally, the mapping $H$ is CR transversal (and hence
transversal) to
$\widetilde M$ at $0$,  by Theorem 1.1 of
\cite {er1}.  The rest of the proof of \eqref{e:equiv2} now
follows from
\eqref{e:equiv1}.

\end{proof}

\begin{Rem} \Label{r:31}  {\rm To prove that $H(M)$ is smooth in the proof above, we have used not only that $\mathcal H(\mathcal M)$ is smooth but also that $\mathcal H$ is transversal to   $\mathcal H(\mathcal M)$. In general, $\mathcal H(\mathcal M)$ smooth does not imply that $H(M)$ is smooth.  For instance, consider $M=\bR$ in $\bC$ and the mapping $z\mapsto z^2$. However, the reverse implication does hold, i.e.\ if $H(M)$ is smooth, then $\mathcal H(\mathcal M)$ is also smooth.
}
\end{Rem}

\section{Further results on images of curves under finite
mappings}\Label{s:last}

In this section, we shall address the following question, which was alluded to above.
\medskip

\noindent {\bf Question:} {\it Let $X$ be a complex submanifold
through $0$ in $\bC^k$, and $f\colon (\bC^k,0)\to (\bC^k,0)$ a germ
of a finite holomorphic mapping. Does the identity $f^{-1}(f(X))=X$,
as germs at 0, imply that $f(X)$ is a submanifold at $0$?}
\medskip

An equivalent formulation can be given as follows:

\noindent {\bf Question$'$:} {\it Let $\tilde X$ be a complex
submanifold through $0$ in $\bC^k$, and $f\colon (\bC^k,0)\to
(\bC^k,0)$ a germ of a finite holomorphic mapping. Assume that
$X:=f^{-1}(\tilde X)$ is a submanifold at $0$. Does this imply that
$\tilde X$ is a submanifold at $0$?}
\medskip

As mentioned above, we do not know the answer in general. However, the answer for
one-dimensional submanifolds is affirmative in view of the following result.

\begin{Thm}\Label{t:curveimage1} Let $X$ be a complex
submanifold of dimension one (i.e. a smooth complex
curve) through $0$ in
$\bC^k$ and $f\colon (\bC^k,0)\to (\bC^k,0)$ a germ of a
finite holomorphic mapping. If $f^{-1}(f(X))=X$ as germs
at $0$, then $\tilde X:=f(X)$ is a germ at $0$ of a submanifold.
\end{Thm}

An equivalent formulation of this result in the spirit of the second
formulation of the question above is the following.

\begin{Thm}\Label{t:curveimage2} Let $\tilde X$ be a
complex subvariety of dimension one (i.e. a complex
curve) through $0$ in
$\bC^k$ and $f\colon (\bC^k,0)\to (\bC^k,0)$ a germ of a
finite holomorphic mapping. If $X:=f^{-1}(\tilde X)$ is
a germ at $0$ of a submanifold (i.e.\ a smooth curve),
then $\tilde X$ is also a submanifold (i.e\ a
smooth curve) at $0$.
\end{Thm}

\begin{proof}[Proof of Theorem $\ref{t:curveimage1}$] Let
$Z=(z,w)\in
\bC\times\bC^{k-1}$ be local coordinates at $0$ in which
$X=\{(z,w)\colon w=0\}$. We also choose coordinates
$\tilde Z$ at $0$ in the target copy of $\bC^k$ such
that, after possibly another change of coordinates in the
$z$ variable,
\begin{equation}
g(z):=f(z,0)=(z^m,f_2(z,0),\ldots, f_k(z,0)),
\end{equation}
where $f_j(z,0)=O(z^m)$ for $j=2,\ldots, k$.  We expand the $f_j(z,0)$, $j=2,\ldots, k$, in their Taylor series
\begin{equation}
f_j(z,0)=\sum_{l=1}^\infty a_{jl} z^l,\quad j=2,\ldots, k.
\end{equation}
We define
\begin{equation}\Label{e:qdef}
q:=\gcd(m, \{l\colon a_{1l}\neq 0\},\ldots, \{l\colon a_{kl}\neq 0\}),
\end{equation}
where $\gcd(n_1,n_2,\ldots)$ denotes the greatest common divisor of the numbers $n_1,n_2,\ldots$.  We observe that there are holomorphic functions $h_j(z)$, for $j=2,\ldots, k$,  such that
\begin{equation}\Label{e:fh}
f_j(z,0)=h_j(z^q),\quad j=2,\ldots, k.
\end{equation}
 We claim that, for any $t\neq 0$, the preimages of $g(t)$  are
 $\{\epsilon_0 t,\ldots, \epsilon_{q-1}t\}$, where $\epsilon_0,\ldots, \epsilon_{q-1}$ are the $q$th roots of unity. The fact that all the points $\epsilon_j t$, $j=0,\ldots, q-1$, are preimages is clear from the definition of $q$ and \eqref{e:fh}. Conversely, if $z_0$ is a preimage of $g(t_0)$ for some $t_0\neq 0$, then $z_0=\epsilon t_0$ for some $m$th root of unity $\epsilon$. Since the only possible preimages of a point $g(t)$ with $t$ close to $t_0$ are of the form $z=\delta t$ for some $m$th root of unity $\delta$, we conclude that  $z=\epsilon t$ is a preimage of $g(t)$ for all $t$ near $t_0$.  Hence, we have $f_j(\epsilon t,0)=f_j(t,0)$ for $j=2,\ldots, k$, which implies
 \begin{equation}
\sum_{l=1}^\infty a_{jl}\epsilon^l t^l=\sum_{l=1}^\infty a_{jl} t^l,\quad j=2,\ldots, k,
\end{equation}
for $t$ near $t_0$. Consequently, $\epsilon^l=1$ for all $l$ such that $a_{jl}\neq 0$ for some $j=2,\ldots, k$. Since $\epsilon ^m=1$ as well, we conclude that $\epsilon ^q=1$, where $q$ is as defined in \eqref{e:qdef}. This proves the claim above.

 If $q=m$, then it is clear
that $\tilde X$ is smooth at $0$. Thus, to complete the
proof of Theorem \ref{t:curveimage1} it suffices to
assume that
$1\leq q<m$ and show that $f^{-1}(f(X))$ contains but is
not equal to $X$. This is an immediate consequence of the
following lemma, since $f(t,0)$ is $O(t^m)$ and $q<m$.
\end{proof}

\begin{Lem}\Label{l:lastlem} Let $f$ and $X$ be as in
Theorem
$\ref{t:curveimage1}$ and $q$ the multiplicity of the
mapping $f|_X\colon X\to \tilde X:=f(X)$.
If $f^{-1}(f(X))=X$ as
germs at $0$, then there is a germ at $0$ of
a holomorphic function $F\colon (\bC^k,0)\to (\bC,0)$
such that $F(f(t,0))=ct^j$, where $c\neq 0$ and $1\leq
j\leq q$. \end{Lem}

\begin{proof} We retain the normalizations as in the
beginning of the proof of Theorem \ref{t:curveimage1}.
Let $p$ be the multiplicity of the finite mapping $f\colon
(\bC^k,0)\to (\bC^k,0)$ and, for generic
$\tilde Z$, let $Z^1:=Z^1(\tilde Z),\ldots,
Z^p:=Z^p(\tilde Z)$ be the preimages of $\tilde Z$ under
$f$. For
$j=1,\ldots, p$, we form the $j$th symmetric combination
of these preimages, i.e.\ $F^j(\tilde Z)=(F^j_1(\tilde
Z),\ldots F^j_k(\tilde Z))$ where
\begin{equation}
F^j_i(\tilde Z)=(-1)^j\sum_{1\leq l_1<\ldots<l_j\leq p}
Z^{l_1}_i\ldots Z^{l_j}_i.
\end{equation}
As is well known, the mappings $F^j$, originally defined
only for generic $\tilde Z$, extend as holomorphic mappings
$(\bC^k,0)\to (\bC^k,0)$. Since all the preimages of
$f(t,0)$, for $(t,0)\in X$, are assumed to lie on $X$
and, hence are of the form $(\epsilon_it,0)$ with
$\epsilon _i^q=1$ for $i=0,1,\ldots, q-1$, we conclude
that
\begin{equation}
F^j(f(t,0))=(c_jt^j,0), \quad j=1,\ldots, p,
\end{equation}
for some constants $c_j$. Thus, to prove Lemma
\ref{l:lastlem}, it suffices to show that $c_j\neq 0$ for
some $j\leq q$. To this end, we introduce the
Weierstrass polynomial
\begin{equation}
P(\tilde Z,x):=x^p+F^1_1(\tilde Z)x^{p-1}+\ldots+
F^{p-1}_1(\tilde Z)x^{1}+ F^p_1(\tilde
Z)=\prod_{l=1}^p(x-Z^l_1),
\end{equation}
where the last identity only holds for generic $\tilde
Z$. We let $R(t,x)$ denote the polynomial $P(f(t,0),x)$
and observe that $R(t,x)$ has the form
\begin{equation}\Label{e:R}
R(t,x)=x^p+c_1tx^{p-1}+\ldots+c_{p-1}t^{p-1}x+c_pt^p.
\end{equation}
Moreover, by construction, the distinct roots of $R(t,x)$
are precisely $x_i=\epsilon_i t$, $i=0,1,\ldots, q$.

\begin{Lem}\Label{l:roots} Let $Q(y)$ be a monic polynomial
of degree
$p$
\begin{equation}\Label{e:Q}
Q(y)=y^p+e_1y^{p-1}+\ldots +e_{p-1}y+e_p.
\end{equation}
If all the roots of $Q(y)$ are $q$-roots of unity with
$q\leq p$, then there is $1\leq j\leq q$ such that
$e_j\neq 0$.
\end{Lem}

\begin{proof} Let $S(w_1,\ldots,w_p,y)$ be the monic polynomial
 in $y$ with polynomial coefficients $c_i(w)$ given by
\begin{equation}
S(w,y) = \prod_{l=1}^p(y-w_i) = y^p+c_1(w)y^{p-1}+\ldots
+c_{p-1}(w)y+c_p(w).
\end{equation}
Note that the polynomial $c_i(w)$ is homogeneous of total degree $i$
in the complex variables $w_1,\ldots,w_p$ and is invariant under
all permutations of the $w_i$.  Furthermore, by the well-known
theorem of elementary invariant theory (see e.g.\ \cite{book}, Theorem 5.3.4), if
$d(w)$ is any polynomial invariant under all permutations of $w$,
then there is a polynomial $a(b_1,\ldots,b_p)$ such that
\begin{equation}\Label{e:invariant}
d(w) = a(c_1(w),\ldots, c_p(w))
\end{equation}
In view of the homogeneity of the $c_i(w)$, it follows from
\eqref{e:invariant} that if $d_q(w)$ is homogeneous of degree $q$,
then there is a polynomial $a_q(b_1,\ldots,b_q)$ (necessarily of
degree $\le q$ and with no constant term) such that
\begin{equation}\Label {e:homog}
d_q(w) = a_q(c_1(w),\ldots, c_q(w))
\end{equation}
Now take $d_q(w):= \sum_{i=1}^p w_i^q$.  If  $w_i= \epsilon_i$ is a
$q$-root of unity for $i = 1,\ldots p$, then
$d_q(\epsilon_1,\ldots,\epsilon_p) = p$, and, in particular, is not
zero.  It follows from \eqref{e:homog} that
$c_{j_0}((\epsilon_1,\ldots,\epsilon_p) \not= 0$ for some $j_0$
with
$1\le j_0 \le q$.  Therefore, since all the roots of $Q(y)$ given by
\eqref{e:Q}  are assumed to be
$q$-roots of unity, it follows that  the coefficient $e_{j_0}$ of
$Q(y)$ is not zero.  This proves Lemma \ref{l:roots}.
\end{proof}

To complete the proof of Lemma \ref{l:lastlem}, we set
$y=x/t$ in \eqref{e:R} and define
$Q(y):=R(t,ty)/t^p$, i.e.\
\begin{equation}
Q(y)=y^p+c_1y^{p-1}+\ldots +c_{p-1}y+c_p.
\end{equation}
All the roots of $Q(y)$ are $q$-roots of unity by
construction. By Lemma \ref{l:roots}, there is $1\leq
j\leq q$ such that $c_j\neq 0$. This proves Lemma
\ref{l:lastlem}.
\end{proof}

We conclude this paper by giving an equivalent algebraic
reformulation of the question posed in the beginning of this
section. Let $I$ be an ideal in $\bC\{Z\}$. Recall that $I$ is the
ideal of a complex analytic subvariety $X$ at $0$ if and only if $I = \sqrt{I}$, i.e $I$ 
is radical. The subvariety $X$ is a submanifold at $0$ if and only
if the ring $\bC\{Z\}/I$ is regular, i.e.\ isomorphic to a power
series ring $\bC\{t\}$.\medskip

\noindent {\bf Question$''$:} {\it Let $\phi\colon \bC\{x_1,\ldots,
x_k\} \to \bC\{z_1,\ldots, z_k\}$ be an injective $\bC$-algebra
homomorphism such that $\bC\{z_1,\ldots, z_k\}$ is integral over
$\phi(\bC\{x_1,\ldots, x_k\})$. Let $I$ be a radical ideal in
$\bC\{x_1,\ldots, x_k\}$ and $J$ the ideal in $\bC\{z_1,\ldots,
z_k\}$ generated by $\phi(I)$. Assume that the ring
$\bC\{z_1,\ldots, z_k\}/\sqrt{J}$ is regular. Does this imply that
$\bC\{x_1,\ldots, x_k\}/I$ is regular? }
\medskip

In this formulation, the answer is negative if the field of complex
numbers is replaced by any field of characteristic $0<p<\infty$ in
view of the following example, communicated to us by Joseph Lipman, who attributed it to Melvin
Hochster. Let $K$ be a field of characteristic
$p$ and consider the homomorphism $\phi\colon K\{u,v,w\}\to
K\{x,y,z\}$ given by $u\mapsto x^p$, $v\mapsto y^p$, and $w\mapsto
z$. If we let $I$ be the (prime) ideal generated by $w^p+uv$, then
$J$ is the ideal generated by $z^p+x^py^p$, which in characteristic
$p$ is equal to $(z+xy)^p$. The radical $\sqrt{J}$ is then generated
by $z+xy$, and $\bC\{z_1,\ldots, z_k\}/\sqrt{J}$ is regular.
However, $\bC\{x_1,\ldots, x_k\}/I$ is not regular. 
 Another counterexample (in characteristic 3) was attributed to
Bill Heinzer. Lipman also informed us that he has an algebraic proof
\cite{Lip} showing that if ``regular" in Question$''$ above is
replaced by ``normal", then the answer is affirmative (indeed, he
proved this statement in a more general context). For an ideal $I$ for which $\dim \bC\{z_1,\ldots, z_k\}/I = 1$, normal is the same as regular and, hence, Lipman's arguments
yield another proof of Theorem \ref{t:curveimage2} above. We would
like to take this opportunity to thank the above mentioned people
for their help and interest in our question.

\section{Acknowledgment}  We would like to thank Dmitri Zaitsev for several  useful remarks concerning an earlier version of this paper.

\def\cprime{$'$}

\end{document}